\documentclass[12pt, letterpaper]{amsart}

\usepackage{amsmath, amssymb}
\usepackage{array}
\usepackage[frame,cmtip,arrow,matrix,line,graph,curve]{xy}
\usepackage{graphpap, color, paralist, pstricks}
\usepackage[mathscr]{eucal}
\usepackage[pdftex]{graphicx}
\usepackage[pdftex,colorlinks,citecolor=blue]{hyperref}
\usepackage{tikz}
\usepackage{tikz-cd}
\usepackage{mathtools}
\usepackage{enumitem}
\usepackage[table]{xcolor}
\usepackage{diagbox}
\usepackage{multirow}

\usepackage[style=alphabetic, url=false, doi=false, isbn=false]{biblatex}
\AtBeginBibliography{\footnotesize}
\makeatletter
\@namedef{subjclassname@2020}{%
  \textup{2020} Mathematics Subject Classification}
\makeatother

\newtheorem{theorem}{Theorem}[section]
\newtheorem{theoremletter}{Theorem}

\newtheorem{proposition}[theorem]{Proposition}
\newtheorem{corollary}[theorem]{Corollary}
\newtheorem{lemma}[theorem]{Lemma}

\newtheorem{conjecture}[theorem]{Conjecture}

\theoremstyle{definition}
\newtheorem{definition}[theorem]{Definition}
\newtheorem{example}[theorem]{Example}

\newtheorem{question}[theorem]{Question}
\newtheorem{remark}[theorem]{Remark}

\newtheorem*{acknowledgement}{Acknowledgement}

\newcommand{\FF}{\mathbb{F} }

\newcommand{\NN}{\mathbb{N} }
\newcommand{\PP}{\mathbb{P} }

\newcommand{\ZZ}{\mathbb{Z} }

\newcommand{\rM}{\mathrm{M} }

\newcommand{\be}{\mathbf{e} }

\newcommand{\bp}{\mathbf{p} }
\newcommand{\bq}{\mathbf{q} }
\newcommand{\br}{\mathbf{r} }

\newcommand{\bX}{\mathbf{X} }

\newcommand{\bx}{\mathbf{x} }
\newcommand{\by}{\mathbf{y} }

\newcommand{\bPhi}{\mathbf{\Phi} }

\def\im{\mathrm{im}\, }

\title{Space-Filling Hypersurfaces}

\author{Victoria Bourgeois-Conti}
\address{Department of Mathematics, Fordham University, New York, NY 10023}

\email{vbb3@fordham.edu}

\author{Han-Bom Moon}
\address{Department of Mathematics, Fordham University, New York, NY 10023}

\email{hmoon8@fordham.edu}

\author{Erick Perez}
\address{Department of Mathematics, Fordham University, New York, NY 10023}

\email{eperez109@fordham.edu}

\author{Amelia Warren}
\address{Department of Mathematics, Fordham University, New York, NY 10023}
\email{aew4@fordham.edu}

\date{\today}

\begin{document}

\begin{abstract}
We seek space-filling hypersurfaces of a projective space over a finite field. They are hypersurfaces contain all rational points. The main goal is to construct a smooth space-filling hypersurface of minimal degree. For each projective space over a finite field, we find a lower bound of the minimum degree of space-filling hypersurface, and construct an example without any singular rational points. 
Finally, over an odd-dimensional projective space or for fields of characteristic two, we provide a construction of a smooth space-filling hypersurface of minimal degree. 
\end{abstract}

\maketitle

\section{Introduction}\label{sec:intro}

Given finite collection of mathematical objects, one may ask two questions from a statistical perspective. What is the average property we may expect? What are the extremal examples? In this paper, we study an extremal problem in algebraic geometry over finite fields. 

Let $\PP^n$ be an $n$-dimensional projective space over a finite field $\FF_q$. If we take a homogeneous polynomial $F \in \FF_q[X_0, X_1, \dots, X_n]$ of degree $d$, it defines a hypersurface $V(F) \subset \PP^n$, and one may count the number of points $|V(F)(\FF_q)|$, where $V(F)(\FF_q)$ is the set of points admitting homogeneous coordinates in $\FF_q$. Such points, the so-called $\FF_q$-rational points, play an important role in the arithmetic of algebraic varieties. Clearly, $0 \le |V(F)(\FF_q)| \le |\PP^n(\FF_q)|$. We say $V(F)$ is \textbf{space-filling} if it attains the maximum $|\PP^n(\FF_q)| = (q^{n+1}-1)/(q-1)$. It is an easy exercise to construct such examples if reducible hypersurfaces are allowed. 

What happens if we impose stronger geometric conditions? Can we find a smooth space-filling hypersurface? The answer given by Poonen in \cite{Poo04} is much stronger: For any smooth projective variety $X$ defined over $\FF_q$, one may find a hypersurface section $Y \subset X$ of sufficiently large degree such that $Y(\FF_q) = X(\FF_q)$. However, his proof is based on a probabilistic argument, hence it does not provide any explicit example. Moreover, it does not tell if one can find a small fixed degree hypersurface. 

\subsection{Main results}

The aim of this paper is to address these questions for $X = \PP^n$. 
\begin{question}
\begin{enumerate}
\item What is the minimum degree of a smooth space-filling hypersurface in $\PP^n$ over $\FF_q$?
\item Can one construct explicit examples attaining this minimum?
\end{enumerate}
\end{question}

The case of odd-dimensional projective space has a very clear criterion. 

\begin{theoremletter}[Corollary \ref{cor:oddnminimaldegree}]\label{thm:mainthmA}
For any odd-dimensional projective space $\PP^n$ over $\FF_q$, the minimum degree of a smooth space-filling hypersurface is $q+1$. 
\end{theoremletter}

We found the statement in \cite[Theorem 1.1]{HK17}, while preparing this manuscript. For the reader's convenience, we include the proof of Theorem \ref{thm:mainthmA} in Section \ref{sec:q+1}.

When $n$ is even, there is no smooth space-filling hypersurface of degree $q+1$ (Corollary \ref{cor:nosmoothq+1}). Finding a general construction working for all $q$ and $n$ seems to be difficult. As a first step, in Section \ref{sec:companion}, we prove the following result. 

\begin{theoremletter}[Theorem \ref{thm:noFqsingpt}]\label{thm:mainthmB}
For any $\PP^n$ with $n \ge 2$, there is a degree $q+2$ space-filling hypersurface without any $\FF_q$-rational singular points. 
\end{theoremletter}

Theorem \ref{thm:mainthmB} and numerical computations in Section \ref{sec:statistics} for small $n$ and $q$ suggest: 

\begin{conjecture}\label{conj:q+2}
For any $n \ge 2$, there is a smooth space-filling hypersurface of degree $q+2$ in $\PP^n$ over $\FF_q$. In particular, the minimum degree of a smooth space-filling hypersurface on an even-dimensional projective space is $q+2$. 
\end{conjecture}

We observe that in odd characteristic, low-dimensional examples may be used to construct examples in higher dimensions (Lemma \ref{lem:dimreduciton}). In Corollary \ref{cor:smhypgen}, we show that for a fixed $\FF_q$, if there are smooth space-filling hypersurfaces of degree $q+2$ for $n = 2, 3, 4$, then there is an example for every $n \ge 2$.

In characteristic two, we prove Conjecture \ref{conj:q+2} in every dimension. 

\begin{theoremletter}[Theorem \ref{thm:char2}]\label{thm:mainthmC}
For any $\FF_q$ of characteristic two and every dimension $n \ge 2$, there is a degree $q+2$ smooth space-filling hypersurface. 
\end{theoremletter}

\subsection{Previous work}

There have been several previous results in this direction. When $n = 2$, Homma and Kim proved that the minimum possible degree is $q+2$, and they classified smooth space-filling curves of this degree \cite{HK13}. In \cite{AG23}, Asgarli and Ghioca constructed smooth space-filling curves of degree $q+3$ for odd prime powers $q < 200$. For odd higher dimensional projective spaces, Theorem \ref{thm:mainthmA} also appears in \cite{HK17, Tir17}. For $\PP^3$, the result was used in a construction of smooth space-filling complete intersection curves \cite{CDMMM23} for finite fields of order $q \le 7$. For a quadric surface $\PP^1 \times \PP^1$, smooth space-filling curves of the smallest bidegree were constructed in \cite{HK23, ACM25}. In $\PP^n$, smooth space-filling hypersurfaces of degree $q+2$ have no transverse $\FF_q$-lines. This property is studied in \cite{ADL21, BCFW25}. 

\subsection{Organization of the paper}

The paper is organized as follows. Section \ref{sec:preliminaries} reviews the background. In Section \ref{sec:q+1}, we show Theorem \ref{thm:mainthmA} and the nonexistence in this degree for even dimensional projective spaces. In Section \ref{sec:q+2}, we fix notations for the investigation of degree $q+2$ space-filling hypersurfaces coming from a square matrix. We argue that constructions in low-dimensional projective space is sufficient in Section \ref{sec:dimreduction}. In Sections \ref{sec:diagonal} and \ref{sec:companion}, we focus on diagonal matrices and companion matrices respectively, and prove Theorem \ref{thm:mainthmB}. Section \ref{sec:statistics} presents computational results. Finally, in Section \ref{sec:chartwo}, we prove Theorem \ref{thm:mainthmC}.

\subsection*{Notation and conventions}
We denote $[n] := \{0, 1, 2, \dots, n\}$. We use uppercase letters for variables and lowercase letters for elements in a field. For any field $k$, $\bar{k}$ denotes its algebraic closure. A fintie field of order $q$ is denoted by $\FF_q$.

\acknowledgement

The second author has been partially supported by the AMS-Simons Research Enhancement Grant for PUI Faculty. The proof of Lemma \ref{lem:geomsmoothchar2} was developed by ChatGPT, and the authors polished the exposition.

\section{Preliminaries}\label{sec:preliminaries}

\subsection{Review of projective varieties}
Let $k$ be a field and $S := k[X_0, X_1, \dots, X_n]$ be the polynomial ring over $k$. The projective space $\PP^n$ is the set of equivalence classes of nonzero  vectors in $\overline{k}^{n+1}$. Two vectors $(x_0, x_1, \dots, x_n)$ and $(y_0, y_1, \dots, y_n)$ are equivalent if they are nonzero scalar multiples of each other, and the equivalence class is denoted by the  homogeneous coordinates $\bx := (x_0 : x_1 : \dots : x_n)$, where $x_i \in \bar{k}$. For an intermediate field $k \subset K \subset \bar{k}$, if we may find a representative $(x_0:x_1:\dots :x_n)$ with $x_i \in K$ for all $0 \le i \le n$, we say $\bx$ is a \textbf{$K$-rational point}. The set of $K$-rational points in $\PP^n$ is denoted by $\PP^n(K)$.

A homogeneous ideal $I \subset S$ defines a \textbf{projective variety} $Y := V(I) \subset \PP^n$, which is the common zero set of polynomials in $I$. To emphasize that the defining equations of $Y$ are polynomials with $k$-coefficients, we say it is a projective variety over $k$. If $I$ is generated by $F_1, F_2, \dots, F_t \in S$, we denote $Y = V(F_1, F_2, \dots, F_t)$. The set of all $K$-rational points of $Y$ is denoted by $Y(K)$ in a similar manner. 

A projective variety that has a single defining equation (hence can be written as $V(F)$) is called a \textbf{hypersurface}. A point $\bx \in Y := V(F)$ is called a \textbf{singular point} if $\nabla F(\bx) = \mathbf{0}$. Otherwise, $\bx$ is called a smooth point. A hypersurface $Y \subset \PP^n$ over $k$ is \textbf{smooth} if every point in $Y(\overline{k})$ is a smooth point. In general, checking smoothness only at $k$-rational points is insufficient. 

We may describe the set of singular points $Y^{\mathrm{sing}}$ as follows. If we denote 
\[
    F_i := \frac{\partial F}{\partial X_i}, 
\]
$Y = V(F)$ is smooth if and only if $Y^{\mathrm{sing}} := V(F, F_0, F_1, \dots, F_n) = \emptyset$. By Euler's identity, 
\[
    \sum_{i=0}^n X_i F_i = (\deg F)F. 
\]
Thus, if $\mathrm{char}\; k$ does not divide $\deg F$, then $F$ is in the ideal generated by its partial derivatives. Hence $Y^{\mathrm{sing}} = V(F_0, F_1, \dots, F_n)$.

\subsection{Space-filling hypersurfaces}\label{ssec:sfh}

In this paper, we are particularly interested in a projective space over a finite field. Let $\FF_q$ be the finite field of order $q$. It is well known that 
\[
    |\PP^n(\FF_q)| = \frac{q^{n+1}-1}{q-1}.
\]
Thus, for any projective variety $X \subset \PP^n$ over $\FF_q$, $|X(\FF_q)| \le (q^{n+1}-1)/(q-1)$. For each dimension $k \le n$, it is easy to construct a subvariety whose number of $\FF_q$-rational points achieves the upper bound. However, if we seek a smooth irreducible variety, this is highly nontrivial. 

A hypersurface $Y = V(F)$ is a \textbf{space-filling hypersurface} if $Y(\FF_q) = \PP^n(\FF_q)$. Our primary interest is the construction of a smooth space-filling hypersurface. Its existence was shown by Poonen, in a more general setting. 

\begin{theorem}[\protect{\cite[Theorem 3.3]{Poo04}}]\label{thm:Poonen}
Let $X \subset \PP^n$ be a smooth geometrically irreducible projective variety of dimension $m \ge 2$ defined over $\FF_q$. Then there exists a smooth hypersurface $H \subset \PP^n$ such that $Y := X \cap H$ is smooth, geometrically irreducible, and $Y(\FF_q) = X(\FF_q)$. 
\end{theorem}

In particular, by applying this argument repeatedly, one may show the existence of a smooth space-filling complete intersection curve $C \subset X$, for any $X$ satisfying these hypotheses, which answers a question of Katz in \cite{Kat99} affirmatively. 

Poonen's argument for Theorem \ref{thm:Poonen} is probabilistic. He showed that the probability that a randomly chosen degree $d$, $\FF_q$-hypersurface $H$ yields a smooth, irreducible intersection $Y := X \cap H$ and $Y(\FF_q) = X(\FF_q)$ is positive, if $d$ is sufficiently large. However, it does not provide any concrete example, nor a lower bound of the degree. 

We now return to the case of $X = \PP^n$. The first non-trivial case is $\PP^2$, where the minimal degree smooth space-filling curves were classified by Homma and Kim in \cite{HK13}. 

\begin{theorem}[\protect{\cite[Theorem 3.2]{HK13}}]\label{thm:HK13}
The minimal degree of a smooth space-filling curve in $\PP^2$ is $q+2$. Any degree $q + 2$ space-filling curve in $\PP^2$ can be written as $V(F)$ with 
\[
    F = \begin{bmatrix}
    X_0& X_1& X_2
    \end{bmatrix}
    A
    \begin{bmatrix}
    X_1^q X_2 - X_1 X_2^q\\
    X_2^q X_0 - X_2 X_0^q\\
    X_0^q X_1 - X_0 X_1^q
    \end{bmatrix}
\]
for some $A \in \rM_{3 \times 3}(\FF_q)$. Moreover, $V(F)$ is smooth if and only if the characteristic polynomial $P_A(T)$ of $A$ is irreducible over $\FF_q$.
\end{theorem}

Our goal is to explore a generalization of Theorem \ref{thm:HK13} to higher-dimensional projective spaces. 

For each dimension, we have a homogeneous ideal generating all space-filling hypersurfaces. Let $S = \FF_q[X_0,X_1, \dots, X_n]$ and $S_d$ be its subspace of homogeneous polynomials of degree $d$. 

\begin{definition}\label{def:phi}
For any two indices $i \ne j$, let
\begin{equation}\label{eqn:phi}
    \phi_{ij} := X_i^q X_j - X_i X_j^q.
\end{equation}
Let $J \subset S$ be the ideal generated by $\phi_{ij}$. We also set $J_k := J \cap S_k$.
\end{definition}

It is straightforward to see that $\phi_{ji} = - \phi_{ij}$. Note also that for $k < q+1$, $J_k = \{0\}$.

\begin{lemma}[\protect{\cite[Lemma 2.4]{CDMMM23}}]\label{lem:idealJ}
A nonzero homogeneous polynomial $F \in S$ defines a space-filling hypersurface if and only if $F \in J$.
\end{lemma}

In particular, if $V(F)$ is a space-filling hypersurface, then $\deg F \ge q+1$.

\section{Hypersurfaces of degree $q + 1$}\label{sec:q+1}

In this section, we give a complete characterization of smooth space-filling hypersurfaces of degree $q+1$. Any degree $q+1$ homogeneous polynomial in $J$ can be written as an $\FF_q$-linear combination 
\begin{equation}\label{eqn:Fdegq+1}
    F = \sum_{0 \le i < j \le n}c_{ij}\phi_{ij}
\end{equation}
for some $c_{ij} \in \FF_q$. 

\begin{definition}\label{def:CF}
For $F$ in \eqref{eqn:Fdegq+1}, we define $C_F \in \rM_{(n+1) \times (n+1)}(\FF_q)$ as
\[
    C_F := 
    \begin{bmatrix}
    0 & -c_{01} & -c_{02} & \dots & -c_{0n}\\
    c_{01} & 0 & -c_{12} & \dots & -c_{1n}\\
    c_{02} & c_{12} & 0 & \dots & \\
    \vdots &&\ddots\\
    c_{0n} & c_{1n} & c_{2n} & \dots & 0
    \end{bmatrix}.
\]
\end{definition}

\begin{lemma}\label{lem:gradient}
The gradient $\nabla F$ is equal to $C_F [X_0^q, X_1^q, \dots, X_n^q]^t$. 
\end{lemma}

\begin{proof}
We compute the partial derivatives. Recall that 
\begin{equation}\label{eqn:derphi}
    \frac{\partial \phi_{ij}}{\partial X_k} = \frac{\partial}{\partial X_k}(X_i^q X_j - X_i X_j^q)
    = \begin{cases}
        -X_j^q, & i = k,\\
        X_i^q, & j = k,\\
        0, & \mbox{otherwise}. 
    \end{cases}
\end{equation}
Then 
\[
    \frac{\partial F}{\partial X_k} = \frac{\partial}{\partial X_k} \sum c_{ij}\phi_{ij} = \frac{\partial}{\partial X_k} (\sum_{i<k} c_{ik} \phi_{ik} + \sum_{k<j} c_{kj}\phi_{kj})
    = \sum_{i < k}c_{ik}X_i^q - \sum_{k < j}c_{kj}X_j^q.
\]
Therefore, we have
\[
    \nabla F = 
    \begin{bmatrix}
    \frac{\partial F}{\partial X_0}\\
    \frac{\partial F}{\partial X_1}\\
    \vdots\\
    \frac{\partial F}{\partial X_n}
    \end{bmatrix}
    = 
    C_F 
    \begin{bmatrix}
    X_0^q\\
    X_1^q\\
    \vdots\\
    X_n^q
    \end{bmatrix}.\qedhere
\]
\end{proof}

\begin{theorem}\label{thm:smoothcrietrionq+1}
A space-filling hypersurface $V(F)$ in $\PP^n$ of degree $q+1$ is smooth if and only if $\det(C_F) \ne 0$. 
\end{theorem}

\begin{proof}
The hypersurface $V(F)$ is singular if and only if there is a point $\bx := (x_0 : x_1 : \dots : x_n) \in V(F)(\overline{\FF}_q)$ such that $\nabla F (\bx) = \mathbf{0}$. By Lemma \ref{lem:gradient}, $C_F\begin{bmatrix}x_0^q, & x_1^q, & \dots, & x_n^q\end{bmatrix}^t = \nabla F(\bx) = \mathbf{0}$ for some point $\bx$. In particular, $C_F$ has a nontrivial kernel, hence $\det (C_F) = 0$.

Conversely, suppose that $\det(C_F) = 0$. Since $C_F$ is a matrix over $\FF_q$, if singular, it has a nonzero $\FF_q$-vector $[x_0,  x_1 ,  \dots ,  x_n]^t$ in the kernel. Note that $x_i^q = x_i$ for any $x_i \in \FF_q$. Thus, $[x_0, x_1, \dots, x_n]^t = [x_0^q, x_1^q, \dots, x_n^q]^t$. If we set $\bx = (x_0 : x_1 : \dots : x_n)$, $\nabla F(\bx) = C_F [x_0^q, x_1^q, \dots, x_n^q]^t = C_F [x_0, x_1, \dots, x_n]^t = \mathbf{0}$. Hence $\bx$ is a singular point of $V(F)$. 
\end{proof}

\begin{remark}
The proof also shows that if $V(F)$ is a singular space-filling hypersurface of degree $q+1$, there is an $\FF_q$-rational singular point.
\end{remark}

The above elementary proof was obtained by us, but we do not claim the originality. While preparing an early draft of this paper, we found Theorem \ref{thm:smoothcrietrionq+1} in \cite[Theorem 1.1]{HK17} and \cite[Proposition 14]{Tir17}.

\begin{example}\label{ex:nonsingalternating}
For any odd $n$, the following alternating matrix
\[
    C_F = \begin{bmatrix}
    0 & -1 & 0 & 0 & \dots & 0 &0\\
    1 & 0 & 0 & 0& \dots & 0&0\\
    0 & 0 & 0 & -1 & \dots & 0&0\\
    0 & 0 & 1 & 0 & \dots & 0& 0\\
    &&&&\ddots &\\
    0 & 0& 0& 0& \dots & 0& -1\\
    0 & 0& 0& 0& \dots & 1 & 0
    \end{bmatrix}
    \in \rM_{(n+1) \times (n+1)}(\FF_q)
\]
has $\det (C_F) = 1$. The associated polynomial 
\[
    F = \phi_{01} + \phi_{23} + \dots + \phi_{n-1,n}
    = X_0^q X_1 - X_0 X_1^q + X_2^q X_3 - X_2 X_3^q + \dots + 
    X_{n-1}^q X_n - X_{n-1} X_n^q
\]
defines a smooth space-filling hypersurface $V(F)$ of degree $q+1$. This example was discussed in \cite{Tal61}. This hypersurface has an interesting geometric property. Indeed, $V(F)$ is not only interpolating all $\FF_q$-rational points in $\PP^n$, but also contains all $\FF_q$-tangent lines on $\PP^n$: any $\FF_q$-line $\ell$ tangent to $V(F)$ is contained in $V(F)$. 
\end{example}

\begin{corollary}\label{cor:oddnminimaldegree}
The minimal degree of a smooth space-filling hypersurface in $\PP^n$ with odd $n$ is $q+1$.
\end{corollary}

\begin{proof}
It is immediate from Lemma \ref{lem:idealJ} and Example \ref{ex:nonsingalternating}.
\end{proof}

\begin{corollary}\label{cor:nosmoothq+1}
Suppose that $n$ is even. Then there is no smooth space-filling hypersurface of degree $q+1$ in $\PP^n$.
\end{corollary}

\begin{proof}
If $n$ is even, for any space-filling hypersurface $V(F)$ of degree $q+1$, $C_F$ is an odd dimensional alternating matrix. Therefore, $\det(C_F)$ is zero. By Theorem \ref{thm:smoothcrietrionq+1}, $V(F)$ is singular. 
\end{proof}

\section{Hypersurfaces of degree $q + 2$}\label{sec:q+2}

From now on, we assume $n \ge 2$. Let $0 \ne F \in J_{q+2}$. It can be written as 
\[
    F = \sum_{i<j}\ell_{ij}\phi_{ij}, 
\]
where $\ell_{ij} = \sum b_{k, ij}X_k$ for some $b_{k,ij} \in \FF_q$. We define a matrix $B \in \rM_{(n+1) \times \binom{n+1}{2}}(\FF_q)$, whose columns are indexed by two-element subsets of $[n]$. Its $(k, ij)$ entry is $b_{k,ij}$. Then 
\[
    F = \bX
    B
    \begin{bmatrix}
    \phi_{01}\\ 
    \phi_{02}\\
    \vdots\\
    \phi_{n-1,n}
    \end{bmatrix},
\]
where $\bX := [X_0, X_1, \dots, X_n]$.

We may simplify the computation by choosing $F \in J_{q+2}$ that can be written as a linear combination of a small number of generators of $J$. Let $I \subset \binom{[n]}{2}$ be an index set. We may take $B_I \in \rM_{(n+1) \times |I|}(\FF_q)$ and denote its $(k,ij)$-entry by $b_{k,ij}$ as before. Then we may define
\begin{equation}\label{eqn:FBI}
    F := \bX B_I
    \begin{bmatrix}
    \phi_{ij}
    \end{bmatrix}_{\{i,j\} \in I.}^t
\end{equation}
To get a smooth hypersurface $V(F)$, $I$ cannot be too small. 

\begin{lemma}\label{lem:singAsingF}
Suppose that $|I| = n+1$ and $B_I \in \rM_{(n+1) \times (n+1)}(\FF_q)$ is singular. Then $V(F)$ has an $\FF_q$-rational singular point. 
\end{lemma}

\begin{proof}
Since $B_I$ is singular, we may choose a nonzero row vector $\by$ such that $\by B_I = \mathbf{0}^t$. 

Let $F_k$ be the $k$-th partial derivative of $F$. By the product rule, 
\begin{equation}\label{eqn:prodrule}
    F_k = \be_k B_I [\phi_{ij}]^t_{\{i,j\}\in I}
    + \bX B_I [\frac{\partial \phi_{ij}}{\partial X_k}]_{\{i,j\}\in I}^t, 
\end{equation}
where $\be_k$ is the $k$-th standard vector. Then 
\[
    F_k (\by) = \be_k B_I [\phi_{ij}(\by)]^t_{\{i,j\}\in I}
    + \by B_I [\frac{\partial \phi_{ij}}{\partial X_k}(\by)]_{\{i,j\}\in I}^t
    = \be_k B_I \mathbf{0} + \mathbf{0}^t
    [\frac{\partial \phi_{ij}}{\partial X_k}(\by)]_{\{i,j\}\in I}^t = 0.
\]
Therefore, $\by \in V(F, F_0, F_1, \dots, F_n)(\FF_q)$. Hence $\by$ defines a singular point of $V(F)$.
\end{proof}

\begin{corollary}\label{cor:smallIsingular}
If $|I| \le n$, then $V(F)$ has an $\FF_q$-rational singular point.
\end{corollary}

\begin{proof}
Suppose that $|I| \le n$. Then by adding $n+1 - |I|$ zero columns to $B_I$ and choosing $n+1-|I|$ new $\phi_{i,j}$'s, we may assume that $F$ is a linear combination of $n+1$ generators. Since we added zero columns on $B_I$, the resulting square matrix is singular. Lemma \ref{lem:singAsingF} gives the result. 
\end{proof}

In particular, to make a smooth hypersurface, we need to use at least $n+1$ generators. One natural choice is $\{\phi_{01}, \phi_{12}, \dots, \phi_{n-1,n}, \phi_{n0}\}$. 

\begin{definition}
We set $\phi_{i} := \phi_{i,i+1}$ with the convention that the indices are taken modulo $n+1$: Thus, $\phi_n = \phi_{n,n+1} = \phi_{n,0}$. 
\end{definition}

Note that if we use exactly $n+1$ generators, the matrix $B_I$ in \eqref{eqn:FBI} is square. For the rest of this paper, we will focus on the following special case. 

\begin{definition}\label{def:matrixnotation}
Let $\bX := [X_0, X_1, \dots, X_n]$ and let $\bPhi := [\phi_{0}, \phi_{1}, \dots, \phi_{n}]$. For any $A \in \rM_{(n+1) \times (n+1)}(\FF_q)$, we define a polynomial 
\begin{equation}\label{eqn:matrixnotation}
    F_A := \bX A\bPhi^t \in J_{q+2}.
\end{equation}
If $F_A \ne 0$, it defines a space-filling hypersurface. We denote the $(i,j)$-th entry of $A$ by $a_{ij}$. The indices are taken modulo $n+1$. 
\end{definition}

Let $\be_i$ be the $i$-th coordinate point of $\PP^n$. It is straightforward to find a condition to have $\be_i$ as a singular point. 

\begin{proposition}\label{prop:coordsing}
The hypersurface $V(F_A)$ in Definition \ref{def:matrixnotation} is singular at $\be_i$ if and only if $a_{i,i-1} = a_{ii} = 0$.
\end{proposition}

\begin{proof}
From the product rule in \eqref{eqn:prodrule} and the partial derivative of $\phi_i$ in \eqref{eqn:derphi}, the $k$-th partial derivative $F_k$ of $F_A$ is 
\[
    F_k = \be_k A [\phi_i]^t + \bX A [\be_{k-1}X_{k-1}^q - \be_k X_{k+1}^q]^t.
\]
So 
\[
    F_k(\be_i) = \begin{cases}
    -a_{i,i-1}, & k = i-1,\\
    a_{ii}, & k = i+1,\\
    0, &\mbox{otherwise}.
    \end{cases}
\]
Therefore, $\nabla F_A(\be_i) = 0$ if and only if $a_{i,i-1} = a_{ii} = 0$.
\end{proof}

Proposition \ref{prop:coordsing} is not used to prove our main theorems. But it motivates our focus on diagonal and companion matrices in the latter sections.

\section{Dimension reduction}\label{sec:dimreduction}

The following observation enables us to reduce the question to lower dimensions. 

\begin{lemma}\label{lem:dimreduciton}
Let $\FF_q$ be a finite field of odd characteristic. Suppose that there are smooth space-filling hypersurfaces $V(F)$ in $\PP^n$ and $V(G)$ in $\PP^m$ of degree $q+2$ over $\FF_q$. Then there is a smooth space-filling hypersurface in $\PP^{n+m+1}$.
\end{lemma}

\begin{proof}
Relabeling the variables, we may assume that $F \in \FF_q[X_0, X_1, \dots, X_n]_{q+2}$ and $G \in \FF_q[Y_0, Y_1, \dots, Y_m]_{q+2}$. Let 
\[
    H := F + G \in \FF_q[X_0, X_1, \dots, X_n, Y_0, Y_1, \dots, Y_m]_{q+2}.
\]

We claim that the hypersurface $V(H) \subset \PP^{n+m+1}$ is smooth and spac-filling. Note that 
\[
    \frac{\partial H}{\partial X_i} = \frac{\partial F}{\partial X_i}, \quad
    \frac{\partial H}{\partial Y_j} = 
    \frac{\partial G}{\partial Y_j}.
\]
Thus, we may write $\nabla H = (\nabla F, \nabla G)$. 

Suppose that $\bp := (q_0 : q_1 : \dots :q_n : r_0 : r_1 : \dots : r_m)\in V(F+G)(\overline{\FF}_q)$ is a singular point, so $\nabla H(\bp) = \mathbf{0}$. Let $\bq := (q_0 : q_1: \cdots : q_n)$ and $\br := (r_0 : r_1 : \cdots : r_m)$. Since 
\[
    \nabla H (\bp) = (\nabla F(\bq), \nabla G(\br)), 
\]
we have $\nabla F(\bq) = \mathbf{0}$ and $\nabla G(\br) = \mathbf{0}$. By Euler's formula, 
\[
    (q+2)F(\bq) = \bq \cdot \nabla F(\bq) = 0, \quad
    (q+2)G(\br) = \br \cdot \nabla G(\br) = 0. 
\]
Since the characteristic is odd, $q+2 \ne 0 \in \FF_q$. Therefore, if $\bq \ne \mathbf{0}$, we have $\bq \in V(F)$. But in this case $\bq$ is a singular point of $V(F)$. By the hypothesis, $V(F)$ is smooth. Thus, $\bq = \mathbf{0}$. Similarly, we have $\br = \mathbf{0}$. Then we have $\bp = \mathbf{0}$, which is a contradiction. 
\end{proof}

Applying Theorem \ref{thm:HK13} and Lemma \ref{lem:dimreduciton} inductively, we obtain the following corollary. 

\begin{corollary}\label{cor:smhyp3n+2}
For any $\FF_q$ of odd characteristic and any $k \in \ZZ_{\ge 0}$, there is a smooth space-filling hypersurface $V(F) \subset \PP^{3k+2}$ of degree $q+2$.
\end{corollary}

Note that the semigroup generated by $3$ and $4$ includes all positive integers $\ge 6$. Thus, we also obtain the following corollary.

\begin{corollary}\label{cor:smhypgen}
Fix a field $\FF_q$ of odd characteristic. 
\begin{enumerate}
\item If there are smooth space-filling hypersurfaces $V(F_k) \subset \PP^k$ of degree $q + 2$ for $k = 2, 3$, then for any $n \ge 5$, there is a smooth space-filling hypersurface $V(F_n) \subset \PP^n$ of degree $q+2$.
\item In addition, if there is a smooth space-filling hypersurface $V(F_4) \subset \PP^4$ of degree $q+2$, then there is a smooth space-filling hypersurface $V(F_n) \subset \PP^n$ of degree $q+2$ for all $n \ge 2$.
\end{enumerate}
\end{corollary}

\section{Diagonal matrices}\label{sec:diagonal}

Among square matrices, perhaps the simplest matrices are diagonal matrices. In this section, we investigate the hypersurfaces $V(F_A)$ of degree $q+2$ associated to a diagonal matrix $A$. The result is negative. Proposition \ref{prop:diagsing} tells us that except for special pairs of $q$ and $n$, all space-filling hypersurfaces of degree $q+2$ from a diagonal matrix have an $\FF_q$-rational singular point. 

Let $A := \mathrm{diag}(a_0,a_1,a_2,\dots, a_n) \in \rM_{(n+1) \times (n+1)}(\FF_q)$. Then 
\[
    F_A = \bX A \Phi^t = 
    \sum_{i=0}^n a_i X_i \phi_i
    = \sum_{i=0}^n a_i X_i (X_i^q X_{i+1} - X_i X_{i+1}^q).
\]

\begin{proposition}\label{prop:diagsing}
Suppose that $q-1$ is relatively prime to $\frac{(-2)^{n+1}-1}{3}$. Then $V(F_A)$ has an $\FF_q$-rational singular point for every diagonal matrix $A$. 
\end{proposition}

Note that $\frac{(-2)^{i}-1}{3} \in \ZZ$ for any $i \in \NN$.

\begin{proof}
If $A$ is singular, $V(F_A)$ has an $\FF_q$-rational singular point by Lemma \ref{lem:singAsingF}. Thus, we may assume that $A$ is invertible, equivalently, $a_i \ne 0$ for all $0 \le i \le n$.

Suppose that $\by = (y_0: y_1:\dots :y_n)$ is an $\FF_q$-rational singular point of $V(F_A)$. We have 
\[
    \frac{\partial}{\partial X_i}F_A
    = a_i \phi_i + a_{i-1}X_{i-1}^{q+1} - a_i X_i X_{i+1}^q.
\]
Since $\phi_i(\by) = 0$ and $y_i^q = y_i$ for $y_i \in \FF_q$, 
\[
    \frac{\partial}{\partial X_i}F_A (\by)
    = a_{i-1}y_{i-1}^{q+1} - a_i y_i y_{i+1}^q
    = a_{i-1}y_{i-1}^2 - a_i y_i y_{i+1}.
\]
Thus, from $\nabla F_A(\by) = \mathbf{0}$, we have `cyclic' recurrence relations
\begin{equation}\label{eqn:singularityVFAdiagonal}
    a_i y_i^2 = a_{i+1} y_{i+1} y_{i+2}, \quad 0 \le i \le n. 
\end{equation}
Here all indices are taken as modulo $n+1$. 

We first claim that this is only possible when every $y_i$ is nonzero. Suppose not, and assume that $y_i = 0$. Then from $a_{i-1} y_{i-1}^2 = a_i y_i y_{i+1} = 0$ and $a_{i-1} \ne 0$, we have $y_{i-1} = 0$. Repeating this argument, we have $y_j = 0$ for all $j$. Hence $\by = \mathbf{0}$, but this is not possible as $\by$ represents a point in $\PP^n$.

Therefore, we may rewrite \eqref{eqn:singularityVFAdiagonal} as: 
\begin{equation}\label{eqn:quotientoftwopriorx}
    y_{i}= \frac{a_{i-2} y_{i-2}^2}{a_{i-1} y_{i-1}}, \quad 0 \le i \le n. 
\end{equation}

Normalizing the vector $\by$, we may assume that $y_0 = 1$ and set $y_1 = t \in \FF_q^*$. Then from \eqref{eqn:quotientoftwopriorx}, we may describe each $y_i$ as a Laurent monomial in $a_0, a_1, \dots, a_i$ and $t$. 

We claim that 
\begin{equation}\label{eqn:yiformula}
    y_i = t^{\frac{1-(-2)^i}{3}}a_0^{\frac{1-(-2)^{i-1}}{3}}\prod_{j=1}^{i-1}a_j^{-(-2)^{i-j-1}}.
\end{equation}
This can be shown by determining the exponent of each variable separately. For instance, if we denote the exponent of $t$ in $y_i$ by $e_i$, from \eqref{eqn:quotientoftwopriorx}, we may obtain
\[
    e_i + e_{i-1} - 2e_{i-2} = 0, \quad e_0 = 0, \quad e_1 = 1.
\]
Solving this recurrence relation, we deduce $e_i = (1-(-2)^i)/3$.

On the other hand, we have $y_{n+1} = y_0 = 1$. Therefore, from \eqref{eqn:yiformula}, we have 
\[
    t^{\frac{(-2)^{n+1}-1}{3}} = a_0^{\frac{1-(-2)^{n}}{3}}\prod_{j=1}^{n}a_j^{-(-2)^{n-j}} =: \alpha.
\]
Thus, $V(F_A)$ has an $\FF_q$-rational singular point $\by$ if $\alpha\in \FF_q^*$ is in the image of the power map $m : \FF_q^* \to \FF_q^*$ given by 
\[
    m(t) = t^{\frac{(-2)^{n+1}-1}{3}}.
\]

Since $\FF_q^*$ is a cyclic group of order $q-1$, $m$ is a group isomorphism if and only if $(q-1, \frac{(-2)^{n+1}-1}{3}) = 1$. Therefore, for any $\alpha \in \FF_q^*$, we have $t \in \FF_q^*$ such that $m(t) = \alpha$. Hence, $V(F_A)$ has an $\FF_q$-rational singular point. 
\end{proof}

\begin{remark}\label{rem:diagexception}
If $(q-1, \frac{(-2)^{n+1}-1}{3}) \ne 1$, the last part of the proof suggests that one may construct a diagonal matrix $A$ such that $V(F_A)$ does not have any $\FF_q$-rational singular points. Moreover, some of them are even smooth. For example, if $n = 3$ and $q = 11$, $(q-1, \frac{(-2)^{n+1}-1}{3}) = 5$. Among diagonal matrices $A = \mathrm{diag}(1,1,1,\lambda)$ for $\lambda \in \FF_{11}$, $\lambda = \pm 2, 3, \pm 4, \pm 5$ yield smooth $V(F_A)$.
\end{remark}

\section{Companion matrices}\label{sec:companion}

The goal of this section is to construct an explicit space-filling hypersurface $V(F)$ of degree $q+2$ for any dimension $n$, without $\FF_q$-rational singular points. Among non-diagonal matrices, another important class consists of companion matrices. A \textbf{companion matrix} is a square matrix of the form 
\begin{equation}\label{eqn:companion}
A = 
\begin{bmatrix}
    0&0& 0&\dots&0&b_0\\
    1&0&0&\dots& 0& b_1\\
    0&1&0&\ddots&0&b_2\\
    \vdots&\vdots&\ddots&\ddots&\vdots&\vdots\\
    0&0&\dots& 1&0&b_{n-1}\\
    0&0&\dots&0&1&b_n
\end{bmatrix}.
\end{equation}
The characteristic polynomial of the companion matrix is
\[
    P_A(T) = T^{n+1} - b_n T^n - \dots - b_1 T - b_0.
\]
In particular, for any monic polynomial $f(T) \in \FF_q[T]$ of degree $n+1$, we may construct a companion matrix $A$ such that $P_A(T) = f(T)$. 

\begin{lemma}\label{lem:singptnonzeroccoeff}
Let $A$ be an invertible companion matrix in \eqref{eqn:companion}. If there is an $\FF_q$-rational singular point $\by = (y_0 : y_1 : \dots : y_n) \in V(F_A)(\FF_q)$, then $y_i \ne 0$ for all $0 \le i \le n$.
\end{lemma}

\begin{proof}
Since $A$ is invertible, $b_0 \ne 0$.  

Suppose that $\by$ is an $\FF_q$-rational singular point. Then $\nabla F_A(\by) = \mathbf{0}$. Using $\phi_i (\by) = 0$ and $y_i^q = y_i$, we obtain the following set of equations:
\begin{equation}\label{eqn:compgrad}
\begin{split}
    y_{n-1} y_n &= y_0 (\sum_{i=0}^n b_i y_i)\\
    y_1^2 &= y_n (\sum_{i=0}^n b_i y_i) \\
    y_k y_{k+1} &= y_{k+2}^2, \quad 0 \le k \le  n-2.
\end{split}
\end{equation}

Suppose $y_0 = 0$. Then from the last equation in \eqref{eqn:compgrad}, $y_2 =  y_3 = \dots = y_n =  0$. The second equation implies $y_1 = 0$, too. Since $\by$ is a point in a projective space, this is not possible. Hence $y_0 \ne 0$. Now assume that $y_1 = 0$. Then using the last equation again, we obtain $y_2 = y_3 = \dots = y_n = 0$. Then the first equation is reduces to $0 = b_0 y_0^2$. Since $b_0 \ne 0$, we have $y_0 = 0$, which is a contradiction. Therefore, $y_0 \ne 0$ and $y_1 \ne 0$. Now from the last equation, we obtain the rest of the coordinates are also nonzero. 
\end{proof}

\begin{proposition}\label{prop:compFqsingpt}
Let $A$ be an invertible companion matrix in \eqref{eqn:companion}. Let $\by$ be an $\FF_q$-rational singular point of $V(F_A)$. Choose homogeneous coordinates of $\by = (y_0:y_1:\dots : y_n)$ with $y_n = 1$. Then there is $z \in \FF_q^*$ such that 
\begin{equation}\label{eqn:compsingcoordy}
    y_k = 
    z^{(1-(-2)^{n-k})/3}
\end{equation}
for every $0 \le k \le n$.
\end{proposition}

\begin{proof}
By Lemma \ref{lem:singptnonzeroccoeff}, $y_i \ne 0$ for all $k$. The last recursive equation in \eqref{eqn:compgrad} can be written as $y_k/y_{k+1} = (y_{k+2}/y_{k+1})^2$. By setting $z_k := y_{k}/y_{k+1}$, we have another form of the recursion
\begin{equation}\label{eqn:compsingrecz}
    z_k = z_{k+1}^{-2}. 
\end{equation}
If we set $z := z_{n-1} = y_{n-1}/y_n = y_{n-1}\in \FF_q^*$, iterating \eqref{eqn:compsingrecz}, we obtain 
\[
    z_k = z^{(-2)^{n-k-1}}
\]
for $0 \le k \le n-1$. Since $y_k = z_k y_{k+1}$ and $y_n = 1$,
\[
    y_k = \prod_{i = k}^{n-1}z_i
    = \prod_{i = k}^{n-1}z^{(-2)^{n-i-1}}
    = z^{(1-(-2)^{n-k})/3}.
\]
The last equality is obtained from summing a finite geometric series. 
\end{proof}

Among invertible companion matrices, we choose one of the simplest ones. 

\begin{definition}\label{def:simplecomp}
Fix $b \in \FF_q^*$. Let $A_{b} \in \rM_{(n+1) \times (n+1)}(\FF_q)$ be a companion matrix such that $b_0 = b$, $b_n = 1$, and $b_j = 0$ for $1 \le j \le n-1$. 
\end{definition}

Note that if $b \in \FF_q^*$, $A_{b}$ is invertible.

\begin{theorem}\label{thm:noFqsingpt}
Set $A := A_{b}$. There is $b \in \FF_q^*$ such that $V(F_A)$ has no $\FF_q$-rational singular point. 
\end{theorem}

\begin{proof}
We use the same normalization as in the proof of Proposition \ref{prop:compFqsingpt}. For an $\FF_q$-rational singular point $\by = (y_0:y_1:\dots : y_n)$, we assume $y_n=1$ and $y_{n-1}=: z$.

The first two equations in \eqref{eqn:compgrad} are simplified as
\[
y_{n-1}y_n = y_0(by_0 + y_n), \quad y_1^2 = y_n(by_0 + y_n).
\]
After applying $y_n = 1$ and $y_{n-1} = z$ and setting $S:= by_0 + 1 \in \FF_q$, they become
\begin{equation}\label{eqn:compgradsim}
z = y_0 S, \quad y_1^2 = S.
\end{equation}
By Proposition \ref{prop:compFqsingpt}, one may check that 
\[
    bz^{\frac{1-(-2)^n}{3}}+1 = by_0 + 1 = S = z/y_0 = z^{1-\frac{1-(-2)^n}{3}}
    = z^{\frac{2+(-2)^n}{3}}.
\]
Therefore, we have an equality
\[
    b = z^{\frac{2(-2)^n+1}{3}}-z^{\frac{(-2)^n-1}{3}}.
\]

Consider a function $\phi : \FF_q^* \to \FF_q$ defined by 
\[
    \phi(t) = t^{\frac{2(-2)^n +1}{3}} - t^{\frac{(-2)^n-1} {3}}.
\]
To have an $\FF_q$-rational singular point $\by$, $b = \phi(z)$ for some $z \in \FF_q^*$. Since $\phi(1) = 0$, $\FF_q^* \setminus \im \phi$ is nonempty. Therefore, if we choose $b \in \FF_q^*\setminus \im \phi$, this is not possible. So with such a choice of $b\in \FF_q^*$, there is no $\FF_q$-rational singular point.
\end{proof}

\begin{remark}
When $\mathrm{char}\; \FF_q = 3$, the examples we constructed in Theorem \ref{thm:noFqsingpt} are not always smooth. We checked for $3 \le n \le 6$ with $q = 3$, and $n = 3$ with $q = 9, 27$. None of the matrices constructed in Theorem \ref{thm:noFqsingpt} yields a  smooth hypersurface. However, except for $\mathrm{char}\; 3$, among hypersurfaces we constructed in Theorem \ref{thm:noFqsingpt}, we don't have any example of non-smooth hypersurfaces. We suspect that $\mathrm{char}\; 3$ is the only exception.
\end{remark}

\section{Computational search}\label{sec:statistics}

We describe computational experiments, particularly for the companion matrices. For each matrix $A$, we construct the corresponding space-filling hypersurface $V(F_A)$ defined in Definition \ref{def:matrixnotation}. Using Macaulay2 \cite{M2}, we tested if the variety is smooth. These computations do not replace the arguments in the previous sections, but they provide examples and reveal patterns. The source code and computational outputs can be found in \cite{BCMPW26}. 

Table~\ref{tbl:q+2} summarizes cases in which the existence of a smooth space-filling hypersurfaces of degree $q+2$ is known. In the table, a colored cell means a smooth space-filling hypersurface in $\PP^n$ over $\FF_q$ exists, hence Conjecture \ref{conj:q+2} is proved for $(n, q)$. 
\begin{enumerate}
\item The red cells labeled H are results of \cite{HK13}. 
\item The blue cells labeled A are established by \cite{ADL21}. We also found new examples using our approach.
\item The purple cells labeled 7 are computational results based on Theorem \ref{thm:noFqsingpt}. 
\item The purple cells labeled 7' are results based on exhaustive search of all invertible companion matrices. Theorem \ref{thm:noFqsingpt} does not provide any smooth example for $\mathrm{char}\; 3$ with $n \ge 3$.
\item The green cells labeled 9 follow from Theorem \ref{thm:mainthmC} proved in Section \ref{sec:chartwo} in this paper. 
\item The cyan cells labeled 5 are results of Section \ref{sec:dimreduction} of this paper.
\item The yellow cells labeled 6 are examples from diagonal matrices. When the coprime assumption in Proposition \ref{prop:diagsing} is violated, we can construct some smooth examples. See Remark \ref{rem:diagexception}. 
\end{enumerate}

\begin{table}[!ht]
\scriptsize
\setlength{\tabcolsep}{3pt} 
\begin{tabular}{|c|c|c|c|c|c|c|c|c|c|c|c|c|c|c|c|c|c|c|c|c|c|c|c|c|c|c|c|c|c|}
\hline
\backslashbox{$n$}{$q$} & 2 & 3 & 4 & 5 & 7 & 8 & 9 & 11 & 13 & 16 & 17 & 19 & 23 & 25 & 27 & 29 & 31 & 32 & 37 & 41 & 43 & 47 & 49 & 53 & 59 & 61 & 64 & 67 & $\cdots$\\ \hline
2 & \cellcolor{red!25}H
& \cellcolor{red!25}H
& \cellcolor{red!25}H
& \cellcolor{red!25}H
& \cellcolor{red!25}H
& \cellcolor{red!25}H
& \cellcolor{red!25}H
& \cellcolor{red!25}H
& \cellcolor{red!25}H
& \cellcolor{red!25}H
& \cellcolor{red!25}H
& \cellcolor{red!25}H
& \cellcolor{red!25}H
& \cellcolor{red!25}H
& \cellcolor{red!25}H
& \cellcolor{red!25}H
& \cellcolor{red!25}H
& \cellcolor{red!25}H
& \cellcolor{red!25}H
& \cellcolor{red!25}H
& \cellcolor{red!25}H
& \cellcolor{red!25}H
& \cellcolor{red!25}H
& \cellcolor{red!25}H
& \cellcolor{red!25}H
& \cellcolor{red!25}H
& \cellcolor{red!25}H
& \cellcolor{red!25}H
& \cellcolor{red!25}H
\\ \hline
3 & \cellcolor{blue!25}A 
& \cellcolor{blue!25}A 
&\cellcolor{green!40}9 
& \cellcolor{blue!25}A 
& \cellcolor{blue!25}A 
&\cellcolor{green!40}9 
& \cellcolor{purple!40}7'
& \cellcolor{blue!25}A 
& \cellcolor{blue!25}A 
&\cellcolor{green!40}9 
& \cellcolor{blue!25}A 
& \cellcolor{blue!25}A 
& \cellcolor{blue!25}A 
& \cellcolor{purple!40}7
& \cellcolor{purple!40}7'
& \cellcolor{blue!25}A 
& \cellcolor{blue!25}A 
&\cellcolor{green!40}9 
& \cellcolor{purple!40}7
& \cellcolor{purple!40}7
& \cellcolor{purple!40}7
&\cellcolor{purple!40}7
&&&
&\cellcolor{yellow!40}6
&\cellcolor{green!40}9 &&
\\ \hline
4 & \cellcolor{green!40}9 
&\cellcolor{purple!40}7' 
& \cellcolor{green!40}9 
&\cellcolor{purple!40}7
&\cellcolor{purple!40}7
& \cellcolor{green!40}9 
&\cellcolor{purple!40}7' 
&\cellcolor{purple!40}7
&\cellcolor{purple!40}7
& \cellcolor{green!40}9 
&&&&&&&& \cellcolor{green!40}9 
&&&&&&&
&
&\cellcolor{green!40}9 
&&
\\ \hline
5 & \cellcolor{green!40}9 
&\cellcolor{cyan!40}5 
& \cellcolor{green!40}9
&\cellcolor{cyan!40}5
&\cellcolor{cyan!40}5
& \cellcolor{green!40}9
& \cellcolor{cyan!40}5
&\cellcolor{cyan!40}5
&\cellcolor{cyan!40}5
& \cellcolor{green!40}9
&\cellcolor{cyan!40}5
&\cellcolor{cyan!40}5
&\cellcolor{cyan!40}5
&\cellcolor{cyan!40}5
&\cellcolor{cyan!40}5
&\cellcolor{cyan!40}5
&\cellcolor{cyan!40}5
& \cellcolor{green!40}9
&\cellcolor{cyan!40}5
&\cellcolor{cyan!40}5
&\cellcolor{cyan!40}5
&\cellcolor{cyan!40}5
&\cellcolor{cyan!40}5
&\cellcolor{cyan!40}5
&\cellcolor{cyan!40}5
&\cellcolor{cyan!40}5
&\cellcolor{green!40}9
&\cellcolor{cyan!40}5
&\cellcolor{cyan!40}5
\\
\hline
6 & \cellcolor{green!40}9 
&\cellcolor{cyan!40}5 
& \cellcolor{green!40}9
&\cellcolor{cyan!40}5
&\cellcolor{cyan!40}5
& \cellcolor{green!40}9
& \cellcolor{cyan!40}5
&\cellcolor{cyan!40}5
&\cellcolor{cyan!40}5
& \cellcolor{green!40}9
&\cellcolor{cyan!40}5
&\cellcolor{cyan!40}5
&\cellcolor{cyan!40}5
&\cellcolor{cyan!40}5
&\cellcolor{cyan!40}5
&\cellcolor{cyan!40}5
&\cellcolor{cyan!40}5
& \cellcolor{green!40}9
&\cellcolor{cyan!40}5
&\cellcolor{cyan!40}5
&\cellcolor{cyan!40}5
&\cellcolor{cyan!40}5
&&&
&\cellcolor{cyan!40}5
&\cellcolor{green!40}9
&&
\\
\hline
7 & \cellcolor{green!40}9 
&\cellcolor{cyan!40}5 
& \cellcolor{green!40}9
&\cellcolor{cyan!40}5
&\cellcolor{cyan!40}5
& \cellcolor{green!40}9
& \cellcolor{cyan!40}5
&\cellcolor{cyan!40}5
&\cellcolor{cyan!40}5
& \cellcolor{green!40}9
&\cellcolor{cyan!40}5
&\cellcolor{cyan!40}5
&\cellcolor{cyan!40}5
&\cellcolor{cyan!40}5
&\cellcolor{cyan!40}5
&\cellcolor{cyan!40}5
&\cellcolor{cyan!40}5
& \cellcolor{green!40}9
&\cellcolor{cyan!40}5
&\cellcolor{cyan!40}5
&\cellcolor{cyan!40}5
&\cellcolor{cyan!40}5
&&&
&\cellcolor{cyan!40}5
&\cellcolor{green!40}9
&&
\\
\hline
8 & \cellcolor{green!40}9 
&\cellcolor{cyan!40}5 
& \cellcolor{green!40}9
&\cellcolor{cyan!40}5
&\cellcolor{cyan!40}5
& \cellcolor{green!40}9
& \cellcolor{cyan!40}5
&\cellcolor{cyan!40}5
&\cellcolor{cyan!40}5
& \cellcolor{green!40}9
&\cellcolor{cyan!40}5
&\cellcolor{cyan!40}5
&\cellcolor{cyan!40}5
&\cellcolor{cyan!40}5
&\cellcolor{cyan!40}5
&\cellcolor{cyan!40}5
&\cellcolor{cyan!40}5
& \cellcolor{green!40}9
&\cellcolor{cyan!40}5
&\cellcolor{cyan!40}5
&\cellcolor{cyan!40}5
&\cellcolor{cyan!40}5
&\cellcolor{cyan!40}5
&\cellcolor{cyan!40}5
&\cellcolor{cyan!40}5
&\cellcolor{cyan!40}5
&\cellcolor{green!40}9
&\cellcolor{cyan!40}5
&\cellcolor{cyan!40}5\\
\hline
9 & \cellcolor{green!40}9 
&\cellcolor{cyan!40}5 
& \cellcolor{green!40}9
&\cellcolor{cyan!40}5
&\cellcolor{cyan!40}5
& \cellcolor{green!40}9
& \cellcolor{cyan!40}5
&\cellcolor{cyan!40}5
&\cellcolor{cyan!40}5
& \cellcolor{green!40}9
&\cellcolor{cyan!40}5
&\cellcolor{cyan!40}5
&\cellcolor{cyan!40}5
&\cellcolor{cyan!40}5
&\cellcolor{cyan!40}5
&\cellcolor{cyan!40}5
&\cellcolor{cyan!40}5
&\cellcolor{green!40}9
&\cellcolor{cyan!40}5
&\cellcolor{cyan!40}5
&\cellcolor{cyan!40}5
&\cellcolor{cyan!40}5
&&&
&\cellcolor{cyan!40}5
&\cellcolor{green!40}9
&&
\\
\hline
$\vdots$ 
& \cellcolor{green!40}9 
&\cellcolor{cyan!40}5 
&\cellcolor{green!40}9
&\cellcolor{cyan!40}5
&\cellcolor{cyan!40}5
& \cellcolor{green!40}9
&\cellcolor{cyan!40}5
&\cellcolor{cyan!40}5
&\cellcolor{cyan!40}5
& \cellcolor{green!40}9
&\cellcolor{cyan!40}5
&\cellcolor{cyan!40}5
&\cellcolor{cyan!40}5
&\cellcolor{cyan!40}5
&\cellcolor{cyan!40}5
&\cellcolor{cyan!40}5
&\cellcolor{cyan!40}5
&\cellcolor{green!40}9
&\cellcolor{cyan!40}5
&\cellcolor{cyan!40}5
&\cellcolor{cyan!40}5
&\cellcolor{cyan!40}5
&&&
&\cellcolor{cyan!40}5
&\cellcolor{green!40}9
&&
\\
\hline
\end{tabular}
\vspace{5pt}
\caption{Existence of smooth space-filling hypersurface in $\PP^n$ over $\FF_q$ of degree $q+2$}
\label{tbl:q+2}
\end{table}

\section{Characteristic two}\label{sec:chartwo}

When $\FF_q$ is of characteristic two, we may even construct a degree $q+2$ smooth hypersurface $V(F)$ for every $n \ge 2$. In this section, $q$ is a power of two. 

A key inductive step is the following proposition. 

\begin{proposition}\label{prop:inductiveconstruction}
Let $\mathrm{char}\;\FF_q = 2$. For $F \in \FF_q[X_0, X_1, \dots, X_n]_{q+2}$, let 
\[
    E := F(X_0, X_1, \dots, X_{n-1}, 0) \in \FF_q[X_0, X_1, \dots, X_{n-1}]_{q+2}.
\]
Assume that $V(F) \subset \PP^{n}$ and $V(E) \subset \PP^{n-1} = \PP^n \cap V(X_n)$ are both smooth space-filling hypersurfaces. Let
\[
    G := F + (X_n + X_{n+1}) \phi_{n, n+1} \in \FF_q[X_0, X_1, \dots, X_{n+1}]_{q+2}
\]
Then $V(G) \subset \PP^{n+1}$ is also a smooth space-filling hypersurface and $G(X_0, X_1, \dots, X_n, 0) = F$.
\end{proposition}

\begin{proof}
Since $G \in J_{q+2}$, $V(G)$ is a space-filling hypersurface. The last assertion is straightforward. 

For a homogeneous polynomial $H$, recall that $H_k$ is its $k$-th partial derivative. Since $\mathrm{char}\; \FF_q = 2$, a direct calculation shows that 
\[
    G_k = \begin{cases}
    F_k, & k < n\\
    F_n + X_{n+1}(X_n + X_{n+1})^q, & k = n\\
    X_n(X_n + X_{n+1})^q, & k = n+1.
    \end{cases}
\]

Suppose that $\bx := (x_0:x_1:\dots:x_{n+1}) \in V(G)(\overline{\FF}_q)$ is a singular point. Let $\bx' := (x_0: x_1 : \dots : x_n)$ and $y := x_n + x_{n+1}$. 

Suppose $y = 0$. Then $F(\bx') = G(\bx) = 0$, hence either $\bx' = \mathbf{0
}$ or $\bx' \in V(F)(\overline{\FF}_q)$. Moreover, $F_k(\bx') = G_k(\bx) = 0$ for all $0 \le k \le n$. If $\bx' \ne \mathbf{0}$, then $\bx'$ is a singular point of $V(F)$, which is not possible by the smoothness assumption of $V(F)$. Thus, $\bx' = \mathbf{0}$. Then $y = x_n = 0$ implies $x_{n+1} = 0$, hence $\bx = \mathbf{0}$, which is a contradiction. 

So $y \ne 0$. From $G_{n+1}(\bx) = 0$, we have $x_n = 0$. If we set $\bx'' := (x_0 : x_1 : \dots : x_{n-1})$, $\bx'' = \mathbf{0}$ or $\bx'' \in V(E)$. If $\bx'' = \mathbf{0}$, then $\bx' = \mathbf{0}$ so $\bx = \be_{n+1}$. In this case, $G_n(\bx) = G_n(\be_{n+1}) \ne 0$, so we have a contradiction. Suppose $\bx'' \in V(E)$. We also have $E_k(\bx'') = F_k(\bx') = 0$ for all $0 \le k \le n-1$. Therefore, $\bx''$ is a singular point of $V(E)$. This is not possible by the assumption. Therefore, $V(G)$ is smooth. 
\end{proof}

Therefore, it is enough to construct a smooth space-filling surface $V(F) \subset \PP^3$ with a smooth hyperplane section. 

Choose $b \in \FF_q^* \setminus \{t^5 + t^3\;|\; t \in \FF_q\}$. Let $A := A_b$ in Definition \ref{def:simplecomp}, that is, a $4 \times 4$ companion matrix 
\[
    \begin{bmatrix}0&0&0&b\\1&0&0&0\\
    0&1&0&0\\0&0&1&1
    \end{bmatrix}
\]
and set $F := F_A$. Note that 
\[
    F = X_1\phi_{01} + X_2 \phi_{12} + X_3 \phi_{23} + (bX_0 + X_3)\phi_{30}. 
\]

By Theorem \ref{thm:noFqsingpt}, $V(F)$ does not have any $\FF_q$-rational singular point. We may strengthen the conclusion using the following lemma. 

\begin{lemma}\label{lem:geomsmoothchar2}
Let $F = F_A$. Then the space-filling surface $V(F)$ is smooth.
\end{lemma}

\begin{proof}
A direct computation yields
\[
    F_i := \frac{\partial F_A}{\partial X_i} = \begin{cases}
    b X_3 X_0^q + X_1^{q+1} + X_3^{q+1},& i = 0,\\
    X_{0}X_1^q+X_{2}^{q+1}, & i = 1, \\
    X_{1}X_2^q+X_{3}^{q+1}, & i = 2, \\
    b X_0^{q+1} + (X_0 + X_2)X_3^q, & i = 3.\end{cases}
\]
Suppose that $\bx = (x_0 : x_1 : x_2 : x_3) \in V(F)(\overline{\FF}_q)$ is a singular point. 

\textbf{Step 1.} If $x_0 = 0$, from the first three partial derivatives, we have $x_1 = x_2 = x_3 = 0$. Therefore, we may set $x_0 = 1$ and $\bx = (1:s:t:u)$. Here a priori, $s, t, u \in \overline{\FF}_q$. From the partial derivatives again, we have
\begin{equation}\label{eqn:singrel}
    s^q = t^{q+1}, \quad s t^q = u^{q+1}, \quad (1+t)u^q = b. 
\end{equation}
In particular, $s, t, u \ne 0$, and $t \ne 1$ from $\mathrm{char}\; \FF_q = 2$.

\textbf{Step 2.} We claim that $t = t^q$. This is equivalent to $t \in \FF_q$. It implies $\bx \in V(F)(\FF_q)$, which contradicts Theorem \ref{thm:noFqsingpt}. Indeed, $u^q = b/(1+t) \in \FF_q$. Then $u^q =: c \in F_q$. But $0 = u^q - c = u^q - c^q = (u-c)^q$, hence $u = c$. Similarly, from $s^q = t^{q+1} \in \FF_q$, we have $s \in \FF_q$. 

\textbf{Step 3.} It remains to prove that $t = t^q$. Set $v := t^q$, $w := v^q = t^{q^2}$, and $z := w^q = t^{q^3}$. Then 
\[
    u^{q(q+1)} = s^q t^{q^2} = t^{q+1}t^{q^2} = tvw.
\]
Since $b \in \FF_q$,
\[
    b = b^q = (1+t)^q u^{q^2} = (1+t^q)u^{q^2} = (1+v)u^{q^2}.
\]
Combining the two equations, we have
\begin{equation}\label{eqn:b2expression}
    b^2 = (1+t)u^q (1+v)u^{q^2}
    = (1+t)(1+v)u^q u^{q^2} = (1+t)(1+v)tvw. 
\end{equation}

\textbf{Step 4.} On the other hand, we have 
\begin{equation}\label{eqn:Fatsing}
\begin{split}
F(1,s,t,u) &= s(s+s^q) + t(s^q t + st^q) + u(t^q u + tu^q) + (b+u)(u^q + u)\\
&= s^2 + s^{q+1} + s^q t^2 + st^{q+1} + t^q u^2 + tu^{q+1} + bu^q + bu + u^{q+1} + u^2.
\end{split}
\end{equation}
From $s^q = t^{q+1}$, we have $s^{q+1} + st^{q+1} = s(s^q+t^{q+1}) = 0$ and $s^q t^2 = t^{q+3}$. Similarly, from $(1+t)u^q = b$, we have $u^{q+1} + tu^{q+1} + bu = 0$.

So the last formula of \eqref{eqn:Fatsing} becomes
\[
    s^2 + t^{q+3} + (t^q +1)u^2 + bu^q
    = s^2 + t^{q+3} + (t+1)^q u^2 + \frac{b^2}{1+t}, 
\]
and the last quantity is zero. Taking its $q$-th power, we have
\[
    s^{2q} + t^{q^2+ 3q} + (t+1)^{q^2} u^{2q} + \frac{b^2}{(1+t)^q} = 0.
\]
Using $s^{2q} = t^{2q+2} = t^2 v^2$ and $u^{2q} = (b/(t+1))^2$, we have
\[
    t^2 v^2 + wv^3 + (w+1)\frac{b^2}{(1+t)^2} + \frac{b^2}{1+v} = 0.
\]
Using \eqref{eqn:b2expression} and multiplying by $(t+1)/v$, we obtain
\[
    t^2 v(t+1) + wv^2 (t+1) + tw(v+1)(w+1) + tw(t+1)^2 = 0.
\]
Expanding it, we have
\[
    t(v+1)w^2 + (t^3 + tv^2 + tv + v^2)w + t^2 v(t+1) = 0.
\]
So if we set 
\[
    R(A, B, C) := A(B+1)C^2 + (A^3 + AB^2 + AB + B^2)C + A^2 B(A+1), 
\]
we have $R(t, v, w) = 0$.

\textbf{Step 5.} As we are working over $\FF_q$, 
\[
    R(v, w, z) = R(t^q, v^q, w^q) = R(t, v, w)^q = 0.
\]
On the other hand, 
\[
    b^2 = b^{2q} = (1+t^q)(1+v^q)t^q v^q w^q = (1+v)(1+w)vwz.
\]
Comparing with \eqref{eqn:b2expression}, we have  
\[
    z = \frac{(1+t)t}{1+w}.
\]
Substituting it for $z$ and a diligent expansion gives
\begin{equation}\label{eqn:identity}
    t(w+1)R(v, w, z)
    + (t^2 + t + v^2)R(t, v, w)
    = w(t+1)(t+v)^4.
\end{equation}
Hence we obtain $w(t+1)(t+v)^4 = 0$. Since we already have $w \ne 0$ and $t\ne 1$, we get $t = v = t^q$ and we are done.
\end{proof}

So we have a good starting point. Note, however, that every restriction of $V(F)$ to a hyperplane of the form $V(X_i)$ is singular. This is because after setting $X_i = 0$, $F$ is a linear combination of only two $\phi_{ij}$, hence singular by Corollary \ref{cor:smallIsingular}. 

To remedy this issue, we use an $\FF_q$-linear change of coordinates of the hypersurface.

\begin{definition}\label{def:H}
For $\lambda \in \FF_q$, let $H(X_0, X_1, X_2, X_3) := F(X_0, X_1, \lambda X_2 + X_3, X_2)$. 
\end{definition}

Since it is obtained from an $\FF_q$-linear change of coordinates, $V(H) \subset \PP^3$ is a smooth space-filling surface.

Now choosing an appropriate $\lambda$, we may obtain a smooth hyperplane section. 

\begin{lemma}\label{lem:char2section}
Let $c \in \FF_q \setminus \{s^3 + s\;|\; s \in \FF_q\}$ and let $\lambda := (c/b)^{q/2}$. Then $V(H) \cap V(X_3) \subset \PP^2$ is also a smooth space-filling curve. 
\end{lemma}

If we define $\psi : \FF_q \to \FF_q$ by $\psi(s) = s^3 + s$, $\psi(0) = \psi(1)$, so $\psi$ is not surjective, hence one can find such $c$.

\begin{proof}
Since $F = X_1\phi_{01} + X_2 \phi_{12} + X_3 \phi_{23} + (bX_0 + X_3)\phi_{30}$, from $\lambda^q = \lambda$, 
\[
    H(X_0, X_1, X_2, 0) = X_1\phi_{01} + \lambda^2 X_2 \phi_{12} + 
    (bX_0 +X_2)\phi_{20}.
\]
With the choice of basis in Theorem \ref{thm:HK13}, it can be written as 
\[
    H(X_0, X_1, X_2, 0) = 
    \begin{bmatrix}X_0 & X_1 & X_2
    \end{bmatrix}
    \begin{bmatrix}
    0 & b & 0\\
    0 & 0 & 1\\
    \lambda^2 & 1 & 0
    \end{bmatrix}
    \begin{bmatrix}
    X_1^q X_2 - X_1 X_2^q\\
    X_2^q X_0 - X_2 X_0^q\\
    X_0^q X_1 - X_0 X_1^q
    \end{bmatrix}.
\]
The characteristic polynomial of the square matrix is 
\[
    P(T) = T^3 + T + \lambda^2 b.
\]

We have $\lambda^2 = (c/b)^q = c/b$ and $\lambda^2 b = c$. By the choice of $c$, we can conclude that $P(T)$ is a cubic polynomial without any $\FF_q$-root. Therefore, $P(T)$ is irreducible in $\FF_q[T]$. By Theorem \ref{thm:HK13}, $V(H(X_0, X_1, X_2, 0)) \subset \PP^2$ is smooth. 
\end{proof}

Using induction on the dimension $n$, we immediately obtain the following result. 

\begin{theorem}\label{thm:char2}
Let $\FF_q$ be a finite field of characteristic two. For every $n \ge 2$, there is a smooth space-filling hypersurface $V(F) \subset \PP^n$ of degree $q+2$.
\end{theorem}

\begin{proof}
The hypersurface $V(H) \subset \PP^3$ in Definition \ref{def:H} provides an example of a smooth space-filling hypersurface of degree $q+2$. The smoothness was shown in Lemma \ref{lem:geomsmoothchar2}, and the smoothness of its section with an appropriate choice of $\lambda$ was shown in Lemma \ref{lem:char2section}. Proposition \ref{prop:inductiveconstruction} shows the existence of a smooth space-filling hypersurface $V(G)$ in the next dimension. Since both $V(G)$ and its hyperplane section $V(H)$ are smooth, we can inductively construct an example in the next dimension. 
\end{proof}

\footnotesize
\printbibliography

\end{document}